\documentclass[12pt]{amsart}
\usepackage{amsmath, amsthm, amssymb, amscd}
\usepackage{graphics}

\title[ ]
{Numerical Radius Norms \\ on Operator Spaces}

\author[T. Itoh]{Takashi Itoh$^*$}
\address{$^*$Department of Mathematics, Gunma University,
         Gunma 371-8510, Japan}
\email{itoh$@$edu.gunma-u.ac.jp}
\author[M. Nagisa]{Masaru Nagisa$^{**}$}
\address{$^{**}$Department of Mathematics and Informatics, Chiba
University,
         Chiba 263-8522, Japan}
\email{nagisa$@$math.s.chiba-u.ac.jp}

\newtheorem{thm}{Theorem}[section]
\newtheorem{prop}[thm]{Proposition}
\newtheorem{lem}[thm]{Lemma}
\newtheorem{cor}[thm]{Corollary}

\theoremstyle{definition}

\newtheorem{defn}[thm]{Definition}

\newtheorem{example}[thm]{Example}

\theoremstyle{remark}

\date{}

\begin{document}

\begin{abstract}

We introduce a numerical radius operator space $(X, \mathcal{W}_n)$.
The conditions to be a numerical radius operator space are weaker than 
the Ruan's axiom for an operator space $(X, \mathcal{O}_n)$.
Let $w(\cdot)$ be the numerical radius norm on $\mathbb{B}(\mathcal{H})$.
It is shown that if $X$ admits a norm $\mathcal{W}_n(\cdot)$ on the matrix space
$\mathbb{M}_n(X)$ which satisfies the conditions, then there is a complete isometry,
in the sense of the norms $\mathcal{W}_n(\cdot)$ and $w_n(\cdot)$,
from $(X, \mathcal{W}_n)$ into $(\mathbb{B}(\mathcal{H}), w_n)$.
We study the relationship between the operator space $(X, \mathcal{O}_n)$ and
the numerical radius operator space $(X, \mathcal{W}_n)$. The category of
operator spaces can be regarded as a subcategory of numerical radius operator
spaces.

\end{abstract}



\vspace{20pt}

\maketitle

\section{Introduction}\label{intro}

\vspace{20pt}
Let $\mathbb{B}(\mathcal{H})$ be the set of all bounded operators 
on a Hilbert space $\mathcal{H}$, and $\mathcal{H}^n$ the $n$-direct sum of
$\mathcal{H}$. We denote by $\| a \|_n$  the operator norm, 
and $w_n(a)$ the numerical radius norm for $a \in \mathbb{B}(\mathcal{H}^n)$ 
respectively, and identify $\mathbb{B}(\mathcal{H}^n)$ with the $n \times n$ matrix
space $\mathbb{M}_n(\mathbb{B}(\mathcal{H}))$.
 
In \cite{ruan}, Ruan introduced a striking concept of operator spaces.
An (abstract) operator space is a complex linear space $X$ together with a
sequence of norms $\mathcal{O}_n(\cdot )$ on the $n \times n$ matrix space 
$\mathbb{M}_n(X)$ 
for each $n \in \mathbb{N}$, which satisfies the following Ruan's axioms 
OI, O${\text{I{\hspace{-1pt}I}}}$:

\begin{align*}
 & {\text O}{\text I}. & \qquad & \mathcal{O}_{m+n}\left( 
\left[ \begin{array}{cc}
x & 0  \\
0   & y  \\
\end{array} \right] \right) 
= \max \{ \mathcal{O}_m(x), \mathcal{O}_n(y) \},
\\
& {\text O}{\text{I{\hspace{-1pt}I}}}.
 & \qquad &\mathcal{O}_n(\alpha x \beta) \le \|\alpha \| \mathcal{O}_m(x) \|\beta \|
\end{align*}
for all $x \in \mathbb{M}_m(X), y \in \mathbb{M}_n(X)$ and
$\alpha \in M_{n,m}(\mathbb{C}), \beta \in M_{m, n}(\mathbb{C})$.

Ruan proved in \cite{ruan} that if $X$ is an (abstract) operator space, 
then there is a complete
isometry $\varPsi$ from $X$ to $\mathbb{B}(\mathcal{H})$, that is,
$\| [\varPsi(x_{ij})] \|_n = \mathcal{O}_n([x_{ij}])$ 
for all $[x_{ij}] \in \mathbb{M}_n(X), n \in \mathbb{N}$. 

In this paper, we introduce an (abstract) numerical radius operator space.
We call that $X$ is a numerical radius operator space 
if  a complex linear space $X$ admits a
sequence of norms $\mathcal{W}_n(\cdot )$ on the $n \times n$ matrix space 
$\mathbb{M}_n(X)$ 
for each $n \in \mathbb{N}$, which satisfies
a couple of conditions ${\text W}{\text I}, {\text W}{\text{I{\hspace{-1pt}I}}}$, where 
${\text W} {\text I}$ is the same as ${\text O}{\text I}$, however 
${\text W}{\text{I{\hspace{-1pt}I}}}$ is a slightly 
weaker condition than ${\text O}{\text{I{\hspace{-1pt}I}}}$ as follows:

\begin{align*}
 & {\text W}{\text I}. & \qquad & \mathcal{W}_{m+n}\left( 
\left[ \begin{array}{cc}
x & 0  \\
0   & y  \\
\end{array} \right] \right) 
= \max \{ \mathcal{W}_m(x), \mathcal{W}_n(y) \},
\\
& {\text W}{\text{I{\hspace{-1pt}I}}}. & \qquad &\mathcal{W}_n(\alpha x \alpha^*) \le 
\|\alpha \|^2 \mathcal{W}_m(x),
\end{align*}
for all $x \in \mathbb{M}_m(X), y \in \mathbb{M}_n(X)$ and
$\alpha \in M_{n,m}(\mathbb{C})$.

It is clear that a subspace $X \subset \mathbb{B}(\mathcal{H})$ is a (concrete) 
numerical radius operator space with $w_n(\cdot)$.

We first show that if $X$ is a numerical radius operator space,
then there is a complete
isometry $\varPhi$, in the sense of norms 
$w_n([\varPhi(x_{ij})])= \mathcal{W}_n([x_{ij}])$ 
for all $[x_{ij}] \in \mathbb{M}_n(X), n \in \mathbb{N}$,
from $X$ to a concrete numerical radius operator space 
in $\mathbb{B}(\mathcal{H})$.

It is well known that there is an equality between the operator norm
and the numerical radius norm so that
$$ \frac{1}{2} \| x\| = w\left( 
\left[ \begin{array}{cc}
0 & x  \\
0   & 0  \\
\end{array} \right] \right) 
\qquad {\text{for}} \: x \in \mathbb{B}(\mathcal{H}).
$$

We next show that, given a numerical radius operator 
space $X$ with $\mathcal{W}_n$, 
defining $\mathcal{O}_n$ by
$$
\text{(OW)} \qquad \frac{1}{2}\mathcal{O}_n(x) =\mathcal{W}_{2n} \left(
\left[ \begin{array}{cc}
0 & x  \\
0   & 0  \\
\end{array} \right] \right)  \qquad {\text{for}}\  x \in \mathbb{M}_n(X),
$$
$X$ becomes an operator space with $\mathcal{O}_n$.  On the other hand,
given an operator space $X$ with $\mathcal{O}_n$, 
the numerical radius operator space which satisfies the equality (OW) 
is not unique.
More precisely, for every operator space $X$, there exists the maximal (resp. minimal) 
numerical radius norm 
$\mathcal{W}_{\text{max}}$
(resp.$\mathcal{W}_{\text{min}}$) affiliated with $(X, \mathcal{O}_n)$ 
(See  the definition in section 3) 
among all of $\mathcal{W}$'s which satisfy WI, 
${\text W}{\text{I{\hspace{-1pt}I}}}$ and (OW). Moreover it is shown that
$\mathcal{W}_{\text{min}} \le \mathcal{W}_{\text{max}} \le 2 \mathcal{W}_{\text{min}}$,
and there are
uncountably many $\mathcal{W}$'s which satisfy WI, 
${\text W}{\text{I{\hspace{-1pt}I}}}$ and (OW) such that 
$$\mathcal{W}_{\text{min}}(x) \le  \mathcal{W}(x) \le \mathcal{W}_{\text{max}}(x)
\qquad {\text{for all }} \: x \in \mathbb{M}_n(X), n \in \mathbb{N}.
$$

  Let $\mathbb{O}$ be the category of the operator spaces in which the objects are operator spaces and
the morphisms are completely bounded maps, $\mathbb{W}$ the category of
the numerical radius operator spaces in which the objects are numerical radius 
operator spaces and the morphisms are $\mathcal{W}$-completely bounded maps
(See the definition in section 2). 
We finally show that $\mathcal{W}_{\text{min}}$ and $\mathcal{W}_{\text{max}}$ are the strict functors which embed
$\mathbb{O}$ into $\mathbb{W}$.

\vspace{20pt}

\section{Numerical radius operator spaces}

In this secton, we are going to prove a representation theorem 
for abstract numerical radius operator spaces. 

 Given abstract numerical radius operator spaces (or operator spaces)
 $X$, $Y$ and a linear map $\varphi$ from $X$ to $Y$, $\varphi_n$ from
$\mathbb{M}_n(X)$ to $\mathbb{M}_n(Y)$ is defined to be
$$
\varphi_n([x_{ij}]) = [\varphi (x_{ij})]
\qquad {\text{for each}} \  [x_{ij}]  \in  \mathbb{M}_n(X), \ n \in \mathbb{N}.
$$ 
We use a simple notation for the norm of $x=[x_{ij}] \in \mathbb{M}_n(X)$
to be $\mathcal{W}(x)$ (resp. $\mathcal{O}(x)$) 
instead of $\mathcal{W}_n(x)$ 
(resp. $\mathcal{O}_n(x)$), 
and for the norm of  $f \in \mathbb{M}_n(X)^*$ to be
$\mathcal{W}^*(f)= 
\sup \{|f(x)| | x = [x_{ij}] \in \mathbb{M}_n(X), \mathcal{W}(x) \le 1 \}$.
We denote the norm of $\varphi_n$ by  $\mathcal{W}(\varphi_n) =
\sup \{ \mathcal{W}(\varphi_n(x)) | x=[x_{ij}]\in  \mathbb{M}_n(X), \mathcal{W}(x) \le 1\}$
(resp. $\mathcal{O}(\varphi_n)=\sup \{ \mathcal{O}(\varphi_n(x)) | x=[x_{ij}]\in  \mathbb{M}_n(X), 
\mathcal{O}(x) \le 1\}$.
The $\mathcal{W}$-completely bounded norm (resp. completely bounded norm) of $\varphi$ 
is defined to be
$$
{\mathcal{W}(\varphi)}_{cb} = \sup \{\mathcal{W}(\varphi_n) | n \in \mathbb{N} \},
$$

$$
(\text{resp.} \quad {\mathcal{O}(\varphi)}_{cb} = \sup \{\mathcal{O}(\varphi_n) | n \in \mathbb{N} \}).
$$
We say $\varphi$ is $\mathcal{W}$-completely bounded (resp. completely bounded) if 
${\mathcal{W}(\varphi)}_{cb} < \infty $ (resp.$  {\mathcal{O}(\varphi)}_{cb} < \infty$), 
and
$\varphi$ is $\mathcal{W}$-completely contractive (resp. completely contractive) if 
${\mathcal{W}(\varphi)}_{cb} \le 1 $ (resp.$  {\mathcal{O}(\varphi)}_{cb} \le 1 $).
We call $\varphi$ is a $\mathcal{W}$-complete isometry (resp. complete isometry) if 
$\mathcal{W}(\varphi_n(x)) =\mathcal{W}(x)$  
(resp.$\mathcal{O}(\varphi_{n}(x)) =\mathcal{O}(x)$)
for each $x  \in  \mathbb{M}_n(X), \ n \in \mathbb{N}$.

The next is fundamental in numerical radius operator spaces like the Ruan's Theorem \cite{ruan}
in the operator space theory.

\begin{thm}
If $X$ is a numerical radius operator space with $\mathcal{W}_n$, then there exist a Hilbert space 
$\mathcal{H}$,
a concrete numerical radius operator space $ Y \subset \mathbb{B}(\mathcal{H})$ with the
numerical radius norm $w(\cdot)$,  
and a $\mathcal{W}$-complete isometry $\varPhi$ from $(X, \mathcal{W}_n)$ 
onto $(Y, w_n)$.
\end{thm}

To prove this theorem, we use the similar argument and idea as in the proof 
of \cite{effrosruan1}. 
We just follow
each step of the proof in \cite{effrosruan1}, 
however we write it down for the convenience of the reader 
because Theorem 2.1 also implies the Ruan's Theorem 
(See Corollary 2.5).
The conditions WI and W${\text{I{\hspace{-1pt}I}}}$ work in the next Lemma. 

\begin{lem}
Let $X$ be a numerical radius operator space.  
If $f \in \mathbb{M}_n(X)^*$ and $\mathcal{W}^*(f) \le 1$, then
there exists a state $p_0$ on $\mathbb{M}_n(\mathbb{C})$ such that
\begin{align*}
&(1)& \quad |f(\alpha x \alpha^*)| &\le & &p_0(\alpha \alpha^*) \mathcal{W}(x),&\\
&(2)& \quad
|f(\alpha x \beta)| &\le & &
2 p_0(\alpha \alpha^*)^{\frac{1}{2}} p_0(\beta^* \beta)^{\frac{1}{2}}
\mathcal{W}\left(
\left[ \begin{array}{cc}
0 & x  \\
0 & 0  \\
\end{array} \right]
\right)& \\
\end{align*}
$
{\text{for all}}
\ \alpha \in \mathbb{M}_{n,r}(\mathbb{C}), x \in \mathbb{M}_{r}(X),
\beta \in \mathbb{M}_{r,n}(\mathbb{C}), r \in \mathbb{N}.
$

\end{lem}

\begin{proof} 
First, we prove the inequality $(1)$. It is sufficient to show the existence
of a state $p_0$ in the state space $S(\mathbb{M}_n(\mathbb{C}))$ of 
$\mathbb{M}_n(\mathbb{C})$ such that

$$
 {\text{Re}} f(\alpha x \alpha^*) \le p_0(\alpha \alpha^*) \mathcal{W}(x)
\qquad {\text{for all}}\ \alpha \in \mathbb{M}_{n,r}(\mathbb{C}), x \in \mathbb{M}_{r}(X).
$$ 

For $\alpha_i \in \mathbb{M}_{n, r_i}(\mathbb{C})$ and
$x_i \in \mathbb{M}_{r_i}(V)$ with 
$\mathcal{W}(x_i) \le 1 \ (i =1, \dots, k, \ k \in \mathbb{N})$,
we define a real valued finction  
$F_{\{\alpha_1, \dots, \alpha_k, x_1, \dots, x_k \}}(\ \cdot \  )$  on 
$S(\mathbb{M}_n(\mathbb{C}))$
by 
$$
F_{\{\alpha_1, \dots, \alpha_k, x_1, \dots, x_k \}} (p) = 
  \sum_{i=1}^k p(\alpha_i \alpha_i^*) - 
{\text{Re}} f(\alpha_i x_i \alpha_i^*), \ {\text{for}}\ p \in S(\mathbb{M}_n(\mathbb{C})).
$$
Set 
$$\bigtriangleup = 
\{ F_{\{\alpha_1, \dots, \alpha_k, x_1, \dots, x_k \}} \mid
\ \alpha_i \in \mathbb{M}_{n, r_i}(\mathbb{C}), 
x_i \in \mathbb{M}_{r_i}(V), \mathcal{W}(x_i) \le 1,  r_i, k \in \mathbb{N} \}.
$$
It is easy to see that  $\bigtriangleup$ is a cone in the set of all real functions on 
$S(\mathbb{M}_n(\mathbb{C}))$.
Let $\bigtriangledown$ be the open cone of all strictly negative functions on 
$S(\mathbb{M}_n(\mathbb{C}))$. 
For any $\alpha_i \in  \mathbb{M}_{n, r_i}(\mathbb{C}), i=1, \cdots, k$, 
there exists $p_1 
\in S(\mathbb{M}_n(\mathbb{C}))$ such that
$p_1(\sum \alpha_i \alpha_i^*) = \| \sum \alpha_i \alpha_i^* \|$.

Since
\begin{eqnarray*}
\mathcal{W}(\sum \alpha_i x_i \alpha_i^*) 
&=& \mathcal{W}\left([\alpha_1, \cdots, \alpha_k] 
\left[ \begin{array}{ccc}
x_1 & &   \\
  & \ddots &  \\
 & & x_k \\
\end{array} \right] 
[\alpha_1, \cdots, \alpha_k]^* \right)
\\
& \le & \| [\alpha_1, \cdots, \alpha_k] \|^2 
\mathcal{W} \left( \left[ \begin{array}{ccc}
x_1 & &   \\
  & \ddots &  \\
 & & x_k \\
\end{array} \right] \right)\\
&=& \| \sum \alpha_i \alpha_i^* \| \max_i \{ \mathcal{W}(x_i) \} \\
& \le & \| \sum \alpha_i \alpha_i^* \|,\\
\end{eqnarray*}
we have 
\begin{eqnarray*}
F_{\{\alpha_1, \dots, \alpha_k, x_1, \dots, x_k \}}(p_1) 
& = & \sum p_1(\alpha_i \alpha_i^*) - {\text{Re}} \sum f(\alpha_i x_i \alpha_i^*) \\
& \ge & \| \sum \alpha_i \alpha_i^* \| - |f(\sum \alpha_i x_i \alpha_i^*)|\\
&\ge &  \| \sum \alpha_i \alpha_i^* \| - \mathcal{W}(\sum \alpha_i x_i \alpha_i^*)\\
&\ge & 0.\\
\end{eqnarray*}
Thus it turns out $\bigtriangleup \cap \bigtriangledown = \o$.

By the Hahn-Banach Theorem,  there exists a measure $\mu $ on $S(\mathbb{M}_n(\mathbb{C}))$ 
such that
$\mu(\bigtriangleup) \ge 0$ and $\mu (\bigtriangledown) < 0$. So we may assume that
$\mu$ is a probability measure.
Set $p_0 = \int p d\mu(p)$.  Since $F_{\{ \alpha, \frac{x}{\mathcal{W}(x)} \}}
 \in \bigtriangleup $, we obtain
$$
p_0(\alpha \alpha^*) - {\text{Re}} f(\alpha \frac{x}{\mathcal{W}(x)} \alpha^*) 
=  \int F_{\{ \alpha, \frac{x}{\mathcal{W}(x)}\}} (p) d\mu(p) \ge 0.
$$

Next, we prove the inequality $(2)$.
Since 
$$ |f(\alpha x \beta)| = 
\left| f \left( [\alpha, \beta^*]
\left[ \begin{array}{cc}
0 & x  \\
0 & 0  \\
\end{array} \right]
[\alpha, \beta^*]^* \right) \right|,
$$
we have 
$$ 
|f(\alpha x \beta)| \le 
p_0(\alpha \alpha^* + \beta^* \beta) 
\mathcal{W} \left(\left[ \begin{array}{cc}
0 & x  \\
0 & 0  \\
\end{array} \right]\right).
$$
Let $\lambda > 0$ and replace 
$\alpha, \beta$ by $\lambda \alpha, \lambda^{-1}\beta$. Then the equality
$$
\inf_{\lambda > 0}
\{
\lambda^2 p_0(\alpha \alpha^*) + \lambda^{-2} p_0(\beta^*\beta)\}
= 2p_0(\alpha \alpha^*)^{\frac{1}{2}}p_0(\beta^* \beta)^{\frac{1}{2}}
$$
implies the desired inequality $(2)$.

\end{proof}

The next  is known as Smith's Lemma \cite{smith} in case that $X$ is an operator space. 

\begin{lem}
Let $X$ be a numerical radius operator space.  
If $\varphi$ is a linear map from $X$ to $\mathbb{M}_n(\mathbb{C})$, then
$$
\mathcal{W}(\varphi)_{cb} = \mathcal{W}(\varphi_n).
$$ 
\end{lem}

\begin{proof}
We can follow the same argument as in the proof of Smith's.
\end{proof}

\begin{lem}
Let $X$ be a numerical radius operator space.  
If $f \in \mathbb{M}_n(X)^*$ and $\mathcal{W}^*(f) \le 1$, then
there exist a $\mathcal{W}$-complete contraction from
$X$ to $\mathbb{M}_n(\mathbb{C})$  and a unit vector $\xi \in (\mathbb{C}^n)^n$
such that
$$
f(x) = (\varphi_n(x) \xi | \xi) \qquad {\text{for all}} \ x \in \mathbb{M}_n(X).
$$
\end{lem}

\begin{proof}
Let $p_0$  be a state which satisfies the inequalities in Lemma 2.2. 
By the GNS construction for $p_0$, we have a representation $\pi$ of
$\mathbb{M}_n(\mathbb{C})$ on a finite dimensional Hilbert space $\mathcal{H}$
and a cyclic vector $\xi_0 \in \mathcal{H}$ such that
$$
p_0(\alpha) = ( \pi (\alpha) \xi_0 | \xi_0).
$$
For $\alpha = [\alpha_1, \dots, \alpha_n] \in \mathbb{M}_{1,n}(\mathbb{C})$, we set
$\tilde{\alpha} =
\left[ \begin{array}{ccc}
\alpha_1 & \cdots & \alpha_n  \\
 &{\text{\Large{0}}}&  \\
\end{array} \right]
 \in \mathbb{M}_n(\mathbb{C})
$
and denote by $\tilde{\mathbb{M}}_n(\mathbb{C})$ all of the elements in the form $\tilde{\alpha}$.
Let $\mathcal{H}_0 = \pi(\tilde{\mathbb{M}}_n(\mathbb{C}))\xi_0$. 
For a fixed $x \in X$, define a quasilinear form $B_x$ on $\mathcal{H}_0 \times \mathcal{H}_0$ by
$$
B_x(\pi(\tilde{\beta})\xi_0, \pi(\tilde{\alpha})\xi_0) = f(\alpha^* x \beta).
$$
Since
\begin{eqnarray*}
|f(\alpha^* x \beta)| 
& \le & p_0(\alpha^* \alpha)^{\frac{1}{2}}p_0(\beta^* \beta)^{\frac{1}{2}} 2 \mathcal{W}
\left(\left[ \begin{array}{cc}
0 & x  \\
0 & 0  \\
\end{array} \right]\right)\\
& =& \|\pi(\tilde{\alpha})\xi_0 \| \|\pi(\tilde{\beta}) \xi_0 \| 2\mathcal{W}
\left(\left[ \begin{array}{cc}
0 & x  \\
0 & 0  \\
\end{array} \right]\right),\\
\end{eqnarray*}
$B_x(\cdot, \cdot)$ is well-defined, and there exists a bounded operator
$\varphi_0(x) \in \mathbb{B}(\mathcal{H}_0)$ such that
$$
f(\alpha^* x \beta) = (\varphi_0(x) \pi(\tilde{\beta}) \xi_0 |  \pi(\tilde{\alpha}) \xi_0).
$$
Since $\dim \mathcal{H}_0 \le n$, we may assume that $\mathcal{H}_0$ is a
subspace of $\mathbb{C}^n$.  Let $e$ be a projection from $\mathbb{C}^n$ onto $\mathcal{H}_0$.
Set $\varphi(x) = \varphi_0(x)e$ for $x \in X$. Then it turns out that 
$\varphi$ maps from $X$ to $\mathbb{M}_n(\mathbb{C})$   and
$$
f(\alpha^* x \beta) = (\varphi(x) \pi(\tilde{\beta})\xi_0 |  \pi(\tilde{\alpha})\xi_0).
$$
We let $e_j=[0, \dots, \stackrel{(j)}{1}, \dots, 0] 
\in \mathbb{M}_{1,n}(\mathbb{C})$ and 
$\xi=
\left[ \begin{array}{c}
\pi(\tilde{e}_1)\xi_0  \\
\vdots  \\
\pi(\tilde{e}_n)\xi_0 \\
\end{array} \right].
$
Then it is not hard to see that
$$
f(x) = (\varphi_n(x) \xi | \xi) \quad {\text{for}} \ x \in \mathbb{M}_n(X) \quad \|\xi \| \le 1.
$$ 

To prove the $\mathcal{W}$-complete boundedness of $\varphi$, by Lemma 2.3, we let
$x =[x_{ij}] \in \mathbb{M}_n(X)$ and
$\xi_1=
\left[ \begin{array}{c}
\pi(\tilde{\gamma}_1)\xi_0  \\
\vdots  \\
\pi(\tilde{\gamma}_n)\xi_0 \\
\end{array} \right] 
$
with $\| \xi_1 \| \le 1$ where $\gamma_i  \in \mathbb{M}_{1,n}(\mathbb{C})$.
Set
$\gamma = \left[ \begin{array}{c}
\gamma_1 \\
\vdots \\
\gamma_n \\
\end{array} \right] \in \mathbb{M}_n(\mathbb{C}).
$
Then it turns out that 
$$p_0(\gamma^* \gamma)  = 
\sum_i  \| \pi(\tilde{\gamma_i}) \xi_0 \|^2 =\| \xi_1 \|^2 \le 1.$$
Thus we have
\begin{eqnarray*}
| (\varphi_n(x)\xi_1 | \xi_1)| 
& =& | \sum ( \varphi_0(x_{ij}) \pi(\tilde{\gamma_j}) \xi_0 |\pi(\tilde{\gamma_j}) \xi_0 )| \\
& =& | \sum f(\gamma_i^* x_{ij} \gamma_j) | \\
& =& |f(\gamma^* x \gamma) |\\
& \le & p_0(\gamma^* \gamma)\mathcal{W}(x)\\
&=& \mathcal{W}(x).
\end{eqnarray*}

\end{proof}

Now we will prove the Theorem 2.1. We denote by $\mathcal{W}CB(X, Y)$ the set of
all $\mathcal{W}$-completely bounded maps from $X$ to a numerical radius operator space $Y$. 

\vspace{20pt}

\begin{flushleft}{\it{Proof of Theorem 2.1}}
\end{flushleft}

Let $\mathcal{C} 
= \cup_{n \in \mathbb{N}} 
\{ \varphi \in \mathcal{W}CB(X, \mathbb{M}_n(\mathbb{C})) \ | \ \mathcal{W}(\varphi)_{cb} \le 1 \}$
and
$\mathcal{H} = \oplus_{\varphi \in \mathcal{C}} \mathbb{C}^{n(\varphi)}$,
where $n(\varphi)$ is the degree of the range space $\mathbb{M}_{n(\varphi)}(\mathbb{C})$
 of $\varphi$.
Define that
$$
\varPhi : X \ni x \longmapsto 
(\varphi (x))_\varphi \in \oplus_{\varphi \in \mathcal{C}} {\mathbb{M}}_{n(\varphi)}(\mathbb{C}).
$$
Since $\mathcal{W}(\varphi)_{cb} \le 1$, it is clear that $\mathcal{W}(\varPhi)_{cb} \le 1$.
Conversely, given any $x \in \mathbb{M}_n(X)$, there exists 
$f \in \mathbb{M}_n(X)^*$ with $\mathcal{W}^*(f) \le 1$ such that $f(x)= \mathcal{W}(x)$ by 
the Hahn-Banach Theorem.  By Lemma 2.4, we find 
$\varphi \in \mathcal{W}CB(X, \mathbb{M}_n(\mathbb{C}))$ with 
$\mathcal{W}(\varphi)_{cb} \le 1$ and a unit vector $\xi \in (\mathbb{C}^n)^n$ 
such that $f(x) =(\varphi_n(x) \xi | \xi)$.
Thus it turns out $ w(\varphi_n(x)) = \mathcal{W}(x)$.
Hence we obtain that $w(\varPhi_n(x)) \ge w(\varphi_n(x)) = \mathcal{W}(x)$.
This completes the proof.

\vspace{20pt}

\begin{cor} \rm{(Ruan's Theorem \cite{ruan})} \quad
If $X$ is an operator space with $\mathcal{O}_n$, then there exist a Hilbert space 
$\mathcal{H}$,
a concrete operator space $ Y \subset \mathbb{B}(\mathcal{H})$,  
and a complete isometry $\varPsi$ from $(X, \mathcal{O}_n)$ onto $(Y, \| \ \|_n)$.
\end{cor}

\begin{proof}

Since $(X, \mathcal{O}_n)$ is also a numerical radius operator space,
we can find a $W$-complete isometry $\varPhi$ 
from $(X, \mathcal{O}_n)$ into $(B(H), w_n)$ by Theorem 2.1. 
We put $\varPsi(x) = \frac{1}{2} \varPhi(x)$.
Then we have for $x \in \mathbb{M}_n(X)$,
\begin{align*}
   & \| \varPsi_n(x)\|_n \le 2 w_n(\varPsi_n(x)) = w_n(\varPhi_n(x))  \\
 = & \mathcal{O}_n(x) = \mathcal{O}_{2n} \left(\begin{bmatrix} x & 0 \\ 0 & 0 \end{bmatrix}\right)
          = \mathcal{O}_{2n} \left(\begin{bmatrix} 0 & x \\ 0 & 0 \end{bmatrix}
               \begin{bmatrix} 0 & 0 \\ 1 & 0 \end{bmatrix}\right)  \\
 \le & \mathcal{O}_{2n} \left(\begin{bmatrix} 0 & x \\ 0 & 0 \end{bmatrix}\right)
         =  w_{2n} \left(\begin{bmatrix} 0 & \varPhi_n(x) \\ 0 & 0 \end{bmatrix}\right)
         =  2 w_{2n}\left(\begin{bmatrix} 0 & \varPsi_n(x) \\ 0 & 0 \end{bmatrix}\right) \\
 = & \| \varPsi_n(x) \|_n .
\end{align*}

\end{proof}

\begin{cor} \quad
If $X$ is a numerical radius operator space with $\mathcal{W}_n$, 
then there exist an operator space norm $\mathcal{O}_n$ on $X$
and a complete $\& $ $\mathcal{W}$-complete isometry 
$\varPhi$ from $X$ into $\mathbb{B}(\mathcal{H})$.
\end{cor}

\begin{proof}
For given $\mathcal{W}_n$ and $x \in \mathbb{M}_n(X)$, we
 define  $\mathcal{O}_n$ 
to be
$\mathcal{O}_n(x) = 2 \mathcal{W}_{2n}
\left(\left[ \begin{array}{cc}
0 & x  \\
0 & 0  \\
\end{array} \right]\right)$.
By Theorem 2.1, there exist a $\mathcal{W}$-complete isometry 
$\varPhi$ from $(X, \mathcal{W}_n)$ 
into $(\mathbb{B}(\mathcal{H}), w_n)$.
Since
$$
\| \varPhi_n(x)\|_n  =2{w}_{2n}
\left(\left[ \begin{array}{cc}
0 & \varPhi_{n}(x) \\
0 & 0  \\
\end{array} \right]\right)
=
2 \mathcal{W}_{2n}
\left(\left[ \begin{array}{cc}
0 & x  \\
0 & 0  \\
\end{array} \right]\right)
= \mathcal{O}_n(x),
$$  $\varPhi$ is also a complete isometry
from $(X, \mathcal{O}_n)$ 
into $(\mathbb{B}(\mathcal{H}), \| \ \|_n)$.

\end{proof}

As in the case of the operator space theory, we can see the basic operations are closed in numerical
radius operator spaces $X, Y$.
For $\varphi =[\varphi_{ij}] \in \mathbb{M}_n(\mathcal{W}CB(X, Y))$, we use the identification
$\mathbb{M}_n(\mathcal{W}CB(X, Y)) =\mathcal{W}CB(X, \mathbb{M}_n(Y))$  by
$\varphi(x) = [\varphi_{ij}(x)]$ for $x \in X$ with the norm $\mathcal{W}(\varphi)_{cb}$.
Especially, $\mathbb{M}_n(X^*)$ is identified with $\mathcal{W}CB(X, \mathbb{M}_n(\mathbb{C}))$
where we give the numerical radius norm $w(\cdot)$ on   $\mathbb{M}_n(\mathbb{C})$.

If $N$ is a closed subspace of $X$, we use the identification 
$\mathbb{M}_n(X/N) = \mathbb{M}_n(X)/\mathbb{M}_n(N)$.

Here we state only the fundamental operations. 

\begin{prop}
Suppose that $X$ and $Y$ are numerical radius operator spaces.
Then 
\begin{enumerate}
\item[(1)] \quad $\mathcal{W}CB(X, Y)$ is a numerical radius operator space.
\item[(2)] \quad The canonical inclusion $X \hookrightarrow X^{**}$ is $\mathcal{W}$-completely
isometric.
\item[(3)] \quad If $N$ is a closed subspace of $X$, then $X/N$ is a numerical radius
operator space.
\end{enumerate}
\end{prop}
\begin{proof}

For (1) and (3),
it is not hard to verify that the norms defined as above on 
$\mathbb{M}_n(\mathcal{W}CB(X, Y))$ and $\mathbb{M}_n(X/N)$
satisfy the conditions
${\text W}{\text{I}}$ and ${\text W}{\text{I{\hspace{-1pt}I}}}$.

To show (2),   since the inclusion 
$i : \mathbb{M}_n(X) \ni x \longmapsto i(x) \in \mathbb{M}_n(X^{**})$
is defined by 
$$<i(x), f> = <f, x> =w([f_{ij}(x_{kl})])\qquad {
\text{for}}\  x \in \mathbb{M}_n(X),\ f \in \mathbb{M}_n(X^*),
$$
we have 
$$
\mathcal{W}(i(x))_{cb} = 
\sup \{ |<f, x>| \  \mid f \in \mathbb{M}_n(X^*), \ \mathcal{W}_{cb}(f) \le 1 \}
= \mathcal{W}(x)
$$
by Lemma 2.4.
\end{proof}

\vspace{20pt}


\section{Numerical radius norms and operator spaces}


In this section, we study the relationship between numerical radius operator
spaces and operator spaces. 

Let $X$ be a numerical radius operator space with $\mathcal{W}_n$. Defining
by 
$\mathcal{O}_n(x) = 2 \mathcal{W}_{2n}
\left(\left[ \begin{array}{cc}
0 & x  \\
0 & 0  \\
\end{array} \right]\right)$ for $x \in \mathbb{M}_n(X)$, $(X, \mathcal{O}_n)$ is an operator space from
Corollary 2.6.  

On the other hand, we let $X$ be an operator space with $\mathcal{O}_n$. 
We call that a
sequence of norms 
$\mathcal{W}_n$ is a numerical radius norm affiliated with $(X, \mathcal{O}_n)$ 
if $\mathcal{W}_n$ satisfies  
${\text W}{\text{I}}$, ${\text W}{\text{I{\hspace{-1pt}I}}}$ and
$$
\text{(OW)} \qquad \frac{1}{2}\mathcal{O}_n(x) =\mathcal{W}_{2n} \left(
\left[ \begin{array}{cc}
0 & x  \\
0   & 0  \\
\end{array} \right] \right)  \qquad {\text{for}} \ x \in \mathbb{M}_n(X).
$$

We often write $\mathcal{W}$(resp. $\mathcal{O}$) instead of $\mathcal{W}_n$
(resp. $\mathcal{O}_n$).

\begin{defn}
We define a norm $\mathcal{W}_{\text{max}}$ on
an operator space $(X, \mathcal{O}_n)$ by 
$$
\mathcal{W}_{\text{max}}(x) = \inf
\frac{1}{2} \| aa^* + b^*b \| 
\quad \text{for} \ x \in \mathbb{M}_n(X),
$$
where  the infimum is taken over all
$a \in \mathbb{M}_{n,r}(\mathbb{C}), y \in \mathbb{M}_r(X), 
b \in \mathbb{M}_{r,n}(\mathbb{C}),  r \in \mathbb{N}$ such that
$x=ayb$ and $\mathcal{O}(y)=1$.

We call $\mathcal{W}_{\text{max}}$ is the maximal numerical 
radius norm affiliated with $(X, \mathcal{O}_n)$.
We note that $a, y, b$ can be chosen from 
$a \in \mathbb{M}_{n}(\mathbb{C}), y \in \mathbb{M}_n(X), 
b \in \mathbb{M}_{n}(\mathbb{C}),  n \in \mathbb{N}$ in the definition of 
$\mathcal{W}_{\text{max}}$
by using the right polar decomposition
of $a= |a^*|u$ and the left polar decomposition of $b = v|b|$. 

\end{defn}

It is easy to see that, for $\ x \in \mathbb{M}_n(X)$, we have
$$
\mathcal{O}(x) = \inf \|a\| \|b\| 
$$
where  the infimum is taken over all
$x=ayb$ as in Definition 3.1.
Then it follows  that
$$
\frac{1}{2}\mathcal{O}(x) \le \mathcal{W}_{\text{max}}(x) \le \mathcal{O}(x) \quad
\text{for} \ x \in \mathbb{M}_n(X).
$$

\begin{thm}
Suppose that $X$ is an operator space with $\mathcal{O}_n$. Then 
 $\mathcal{W}_{\text{max}}$ is a numerical radius norm affiliated with $(X, \mathcal{O}_n)$ 
and the maximal 
among all of numerical radius norms affiliated with $(X, \mathcal{O}_n)$.
\end{thm}

\begin{proof}
First we show that $\mathcal{W}_{\text{max}}$ is a norm.
To see that $\mathcal{W}_{\text{max}}(x_1 + x_2) \le 
\mathcal{W}_{\text{max}}(x_1) + \mathcal{W}_{\text{max}}(x_2)$ 
for $x_1, x_2 \in \mathbb{M}_n(X)$, let
$x_i = a_i y_i b_i, \mathcal{O}(y_i) =1 (i =1,2)$.
Since
$$
x_1 + x_2 = [a_1, a_2] 
\left[\begin{array}{cc}
y_1 & 0 \\
0 & y_2 \\
\end{array}\right]
\left[\begin{array}{c}
b_1\\
b_2\\
\end{array}\right]  \ {\text{and}} \ 
\mathcal{O}\left(
\left[\begin{array}{cc}
y_1 & 0 \\
0 & y_2 \\
\end{array}\right]\right) = 1,
$$
we have
\begin{align*}
\mathcal{W}_{\text{max}}(x_1 + x_2) & \le 
\frac{1}{2} \|a_1a_1^* + a_2a_2^* + b_1^*b_1 + b_2^*b_2 \|\\
& \le 
\frac{1}{2} \|a_1a_1^* + b_1^*b_1\| +\frac{1}{2} \|a_2a_2^* + b_2^*b_2 \| .\\
\end{align*}
It is easy to show the rest of the norm conditions.

Next we prove that $ \mathcal{W}_{\text{max}}$ satisfies ${\text W}{\text{I}}$
and ${\text W}{\text{I{\hspace{-1pt}I}}}$.
To see ${\text W}{\text{I}}$, let 
$\left[\begin{array}{cc}
x_1 & 0 \\
0 & x_2 \\
\end{array}\right]
= ayb
$ and $\mathcal{O}(y) = 1$. Since
$x_1 = [1, 0] ayb \left[\begin{array}{c}1\\ 0\\ \end{array}\right]$, we have
\begin{align*}\mathcal{W}_{\text{max}}(x_1) & 
\le \frac{1}{2} 
\|
[1, 0] aa^* \left[\begin{array}{c}1\\ 0\\ \end{array}\right] +
[1, 0] b^*b \left[\begin{array}{c}1\\ 0\\ \end{array}\right]
\| \\
&
\le \frac{1}{2}\|aa^* + b^*b\|. \\
\end{align*}
Also we have $\mathcal{W}_{\text{max}}(x_2) \le \frac{1}{2}\|aa^* + b^*b\|$.
Thus it turns out that
$$ \mathcal{W}_{\text{max}}\left( 
\left[ \begin{array}{cc}
x_1 & 0  \\
0   & x_2  \\
\end{array} \right] \right) 
\le \max \{ \mathcal{W}_{\text{max}}(x_1), \mathcal{W}_{\text{max}}(x_2) \}.
$$
Conversely, let $x_i = a_i y_i b_i, \mathcal{O}(y_i) =1 \ (i=1, 2)$.
Since
$$
\left[\begin{array}{cc}
x_1 & 0 \\
0 & x_2 \\
\end{array}\right] =
\left[\begin{array}{cc}
a_1 & 0 \\
0 & a_2 \\
\end{array}\right]
\left[\begin{array}{cc}
y_1 & 0 \\
0 & y_2 \\
\end{array}\right]
\left[\begin{array}{cc}
b_1 & 0 \\
0 & b_2 \\
\end{array}\right],
$$
we have
\begin{align*} \mathcal{W}_{\text{max}}\left( 
\left[ \begin{array}{cc}
x_1 & 0  \\
0   & x_2  \\
\end{array} \right] \right) & \le
\frac{1}{2} \left\|
\left[\begin{array}{cc}
a_1a_1^* + b_1^*b_1  & 0 \\
0 & a_2a_2^* + b_2^*b_2 \\
\end{array}\right] \right\| \\
& \le 
\max \{ \frac{1}{2} \|a_1a_1^* + b_1^*b_1 \|, \frac{1}{2} \|a_2a_2^* + b_2^*b_2 \| \}.\\
\end{align*}

To see ${\text W}{\text{I{\hspace{-1pt}I}}}$, let 
$x=ayb, \mathcal{O}(y) =1$ and $\alpha \in \mathbb{M}_n(\mathbb{C})$.
Then
\begin{align*}
\mathcal{W}_{\text{max}}(\alpha x \alpha^*) &
\le \frac{1}{2} \| \alpha a a^* \alpha^* + \alpha b^*b \alpha^* \| \\
&
= \frac{1}{2} \| \alpha (a a^* + b^*b) \alpha^* \| \\
&
\le \frac{1}{2} \| a a^* + b^*b \| \|\alpha \alpha^* \| .\\
\end{align*}

To see the condition (OW), let $\mathcal{O}(x) =1$.
Since 
$$
\left[\begin{array}{cc}
0 & x \\
0 & 0 \\
\end{array}\right] =
\left[\begin{array}{cc}
1 & 0 \\
0 & 0 \\
\end{array}\right]
\left[\begin{array}{cc}
x & 0 \\
0 & 0 \\
\end{array}\right]
\left[\begin{array}{cc}
0 & 1 \\
0 & 0 \\
\end{array}\right] \quad {\text{and}} \
\mathcal{O}\left(\left[\begin{array}{cc}
x & 0 \\
0 & 0 \\
\end{array}\right]\right) =1,
$$
we have 
$$
\mathcal{W}_{\text{max}}\left( 
\left[ \begin{array}{cc}
0 & x  \\
0 & 0  \\
\end{array} \right] \right) \le
\frac{1}{2} \left\|
\left[\begin{array}{cc}
1 & 0 \\
0 & 0 \\
\end{array}\right] 
\left[\begin{array}{cc}
1 & 0 \\
0 & 0 \\
\end{array}\right]^* +
\left[\begin{array}{cc}
0 & 1 \\
0 & 0 \\
\end{array}\right]^*
\left[\begin{array}{cc}
0 & 1 \\
0 & 0 \\
\end{array}\right] \right\| = \frac{1}{2}.
$$
To get the other  inequality, let $ x \in \mathbb{M}_r(X)$ with $\mathcal{O}(x) =1$. 
By the Ruan's Theorem, there exist a complete isometry 
$\varphi : X \longrightarrow \mathbb{B}( \mathcal{H})$.
Given $\varepsilon > 0$, we find a unit vectors $\xi, \eta \in \mathcal{H}^r$ such that
$ 1 - \varepsilon < ( \varphi_r(x) \xi | \eta).$
Define $F \in \mathbb{M}_{2r}(X)^*$ 
{for} 
$ 
\left[\begin{array}{cc}
x_1 & x_2 \\
x_3 & x_4 \\
\end{array}\right] 
\in \mathbb{M}_{2r}(X)
$
by
$$
F\left(\left[\begin{array}{cc}
x_1 & x_2 \\
x_3 & x_4 \\
\end{array}\right]\right)
=\left(
\left[\begin{array}{cc}
\varphi_r(x_1) & \varphi_r(x_2) \\
\varphi_r(x_3) & \varphi_r(x_4) \\
\end{array}\right] 
\left[\begin{array}{c}
0 \\
\xi \\ 
\end{array}\right] \ 
| \ 
\left[\begin{array}{cc}
\eta \\
0 \\
\end{array}\right] 
\right) .
$$
We show that $\mathcal{W}_{\text{max}}^*(F) \le 2$.  
Let $z \in \mathbb{M}_{2r} (X)$ with $\mathcal{W}_{\text{max}}(z) <1.$
We may assume that $ z = ayb, \mathcal{O}(y) =1$ and 
$\|aa^* + b^*b \| < 2$ where $ y \in \mathbb{M}_k(X), 
a \in \mathbb{M}_{2r,k}(\mathbb{C})$ and $b \in \mathbb{M}_{k, 2r}(\mathbb{C}).$
Since
\begin{align*}
F(z) & = ( a \varphi_k(y) b 
\left[\begin{array}{c}
0 \\
\xi\\
\end{array} \right] |
\left[\begin{array}{c}
\eta\\
0\\
\end{array} \right])\\
& =
\left(
\left[\begin{array}{cc}
0 & \varphi_k(y)\\
0 & 0 \\
\end{array}\right]
\left[\begin{array}{c}
a^*\\
b\\
\end{array}\right]
\left[\begin{array}{c}
0\\
\xi\\
\end{array}\right] | 
\left[\begin{array}{c}
a^*\\
b\\
\end{array}\right] 
\left[\begin{array}{c}
\eta\\
0\\
\end{array}\right]
\right)\\
& \le \| \varphi_k(y) \| 
\left\|
\left[\begin{array}{c}
a^*\\
b\\
\end{array}\right]
\right\|^2 < 2, \\
\end{align*}
we obtain that
\begin{align*}
\mathcal{W}_{\text{max}} \left( \left[
\begin{array}{cc}
0 & x\\
0 & 0\\
\end{array}\right]\right) & \ge 
\frac{1}{2} F 
\left( \left[
\begin{array}{cc}
0 & x\\
0 & 0\\
\end{array}\right]\right) \\
& =
\frac{1}{2}\left|\left(
\left[\begin{array}{cc}
0 & \varphi_r(x) \\
0 & 0 \\
\end{array}\right] 
\left[\begin{array}{c}
0 \\
\xi \\ 
\end{array}\right] \ 
| \ 
\left[\begin{array}{cc}
\eta \\
0 \\
\end{array}\right] 
\right)\right| \\
& =\frac{1}{2} |(\varphi_r(x)\xi | \eta)| \\
& > \frac{1 - \varepsilon}{2}\\
\end{align*}

Finally we show the maximality of $\mathcal{W}_{\text{max}}$ in the set
of all numerical radius norms affiliated with $(X, \mathcal{O}_n)$.  
To see this, let $\mathcal{W}$ be an
arbitrary numerical radius norm affiliated with $(X, \mathcal{O}_n)$ and
$x = a y b, y \in \mathbb{M}_k(X), 
a \in \mathbb{M}_{n,k}(\mathbb{C})$ and $b \in \mathbb{M}_{k, n}(\mathbb{C}).$
Then we have
\begin{align*}
\mathcal{W}(x) & = \mathcal{W}\left([a, b^*]
 \left[
\begin{array}{cc}
0 & y\\
0 & 0\\
\end{array}\right]
\left[
\begin{array}{c}
a^*\\
b\\
\end{array}\right] \right) \\
&
\le \|[a, b^*] \|^2 \mathcal{W}\left(\left[
\begin{array}{cc}
0 & y\\
0 & 0\\
\end{array}\right]\right) \\
&
= \frac{1}{2} \| aa^* + b^*b \| \mathcal{O}(y) \\
\end{align*}
This implies that $\mathcal{W}(x) \le \mathcal{W}_{\text{max}}(x)$ and completes the proof.
\end{proof}

Next we set $\mathcal{W}_{\text{min}}(x) = \frac{1}{2}\mathcal{O}(x)$ for $x \in \mathbb{M}_n(X).$
It is clear that $\mathcal{W}_{\text{min}}$ satisfies
${\text W}{\text{I}}$, ${\text W}{\text{I{\hspace{-1pt}I}}}$ and (OW).
We can characterize numerical radius norms affiliated with an operator
space $X$ by using $\mathcal{W}_{\text{min}}$ and $\mathcal{W}_{\text{max}}$. We call 
$\mathcal{W}_{\text{min}}$ is the minimal numerical radius norm affiliated with 
$(X, \mathcal{O}_n)$.

\begin{cor}
Suppose that $X$ is an operator space with $\mathcal{O}_n$,  and $\mathcal{W}_n$ satisfies
${\text W}{\text{I}}$, ${\text W}{\text{I{\hspace{-1pt}I}}}$. Then the 
following are equivalent:
\begin{enumerate}
\item[(1)] \quad \rm{(OW)} \quad 
$\frac{1}{2}\mathcal{O}_n(x) =\mathcal{W}_{2n} \left(
\left[ \begin{array}{cc}
0 & x  \\
0   & 0  \\
\end{array} \right] \right)  \quad {\text{for}}\  x \in \mathbb{M}_n(X),
$
\item[(2)]\quad 
There exists a complete and $\mathcal{W}$-complete isometry 
$\varPhi : X \longrightarrow \mathbb{B}(\mathcal{H})$, 
\item[(3)]\quad 
$\mathcal{W}_{\text{min}}(x) \le \mathcal{W}(x) \le \mathcal{W}_{\text{max}}(x) \quad 
{\text{for}}\  x \in \mathbb{M}_n(X).$
\end{enumerate}
\end{cor}
\begin{proof}
(1) $\Rightarrow$ (2)
It follows from the same argument as in the proof of Corollary 2.6.

(2) $\Rightarrow$ (3)
Let $x \in \mathbb{M}_n(X)$. Then we have
$
\mathcal{W}_{\text{min}}(x) = \frac{1}{2}\mathcal{O}(x)
= \frac{1}{2} \| \varPhi_n(x)\| \le w(\varPhi_n(x)) = \mathcal{W}(x)
$
and
$
\mathcal{W}(x) \le \mathcal{W}_{\text{max}}(x)
$
by Theorem 3.2.

(3) $\Rightarrow$ (1)
Let $x \in \mathbb{M}_n(X)$. Then we have

$$
\frac{1}{2}\mathcal{O}(x) = 
\mathcal{W}_{\text{min}}
\left( \left[
\begin{array}{cc}
0 &x \\
0 & 0\\
\end{array} \right] \right)
\le
\mathcal{W}\left( \left[
\begin{array}{cc}
0 &x \\
0 & 0\\
\end{array} \right] \right) \le
\mathcal{W}_{\text{max}}
\left( \left[
\begin{array}{cc}
0 & x \\
0 & 0\\
\end{array} \right] \right)
=\frac{1}{2}\mathcal{O}(x).
$$

\end{proof}

\begin{example}
Let $(X, \mathcal{O}_n)$ be an operator space.  
We present that there are uncountably many numerical radius norms
affiliated with $(X, \mathcal{O}_n)$.

From Corollary 3.3, there exists a complete and $\mathcal{W}$-complete
isometry $\varPhi_{\text{max}} : X \longrightarrow \mathbb{B}(\mathcal{H})$ when we introduce 
the maximal numerical radius norm $\mathcal{W}_{\text{max}}$ on $X$.
Let  $0 \le t \le 1$ and 
$$
a_t =\left[
\begin{array}{ccccc}
0 & 1 &    &  &  \\
   & 0 & t &  &  \\
   &    & \ddots &  \ddots& \\
   &    &          &  \ddots & t \\
   &    &           &  & 0 \\
\end{array}
\right]  \in \mathbb{M}_n(\mathbb{C}), \quad n \ge 3. 
$$
Define that $\varPhi_t(x) = \varPhi_{\text{max}}(x) \otimes a_t$ for $x \in X$. Since $\|a_t\| = 1$, 
then
$\varPhi_t : X \longrightarrow \mathbb{B}(\mathcal{H})\otimes \mathbb{M}_n(\mathbb{C})$ is
completely isometric.
Set $\mathcal{W}_t(x) = w_m([\varPhi_t(x_{ij})])$ for $x=[x_{ij}] \in \mathbb{M}_m(X)$.
It is clear that $\mathcal{W}_t$ is a numerical radius norm affiliated with 
$(X, \mathcal{O}_n)$.
We show that 
$$
\mathcal{W}_{\text{max}}(x) \cos{\frac{\pi}{n+1}} \le \mathcal{W}_1(x) 
\le \mathcal{W}_{\text{max}}(x) \quad {\text{for}} \ x \in \mathbb{M}_m(X), \ m \in \mathbb{N}.
$$
To see this,  given $x=[x_{ij}] \in \mathbb{M}_m(X)$ and $\varepsilon > 0$. Then there exists
a unit vector $\xi \in \mathcal{H}^m$ such that
$|([\varPhi_{\text{max}}(x_{ij})] \xi \ | \ \xi) |> w([\varPhi_{\text{max}}(x_{ij})]) - \varepsilon.$
From \cite{haagerupdela}, we can find a unit vector $\eta \in \mathbb{C}^n$ such that
$w(a_1) = |(a_1 \eta \ | \ \eta)| = \cos{\frac{\pi}{n+1}}$.
Then we obtain that
\begin{align*}
\mathcal{W}_1(x) 
& \ge |(([\varPhi_{\text{max}}(x_{ij}) ]\otimes a_1) \xi \otimes \eta \ | \ \xi \otimes \eta )|\\
& = |([\varPhi_{\text{max}}(x_{ij})] \xi \ | \ \xi)| |(a_1 \eta \ | \ \eta)| \\
& > (\mathcal{W}_{\text{max}}(x) - \varepsilon)\cos{\frac{\pi}{n+1}}.\\
\end{align*}
This implies that $\mathcal{W}_{\text{max}}(x) \cos{\frac{\pi}{n+1}} \le \mathcal{W}_1(x) $.
The second inequality is clear because of the maximality of $\mathcal{W}_{\text{max}}$.
We note that $\mathcal{W}_0 =\mathcal{W}_{\text{min}}$.
Since $[0, 1]\ni t \longmapsto \mathcal{W}_t(x) \in \mathbb{C}$ is continuous, 
then there exist
uncountably many distinct numerical radius norms 
$\{\mathcal{W}_t\}_{0 \le t \le 1}$ affiliated with $(X, \mathcal{O}_n)$. 

There are many ways to construct the numerical radius norms like 
$\{\mathcal{W}_t\}_{0 \le t \le 1}$ affiliated with $(X, \mathcal{O}_n)$.
For instance, replace $a_t$ by
$$
b_t =\left[
\begin{array}{cc}
0 & \sqrt{1 - t} \\
0  &  \sqrt{t}\\
\end{array}
\right]  \in \mathbb{M}_2(\mathbb{C}).
$$

\end{example}

\begin{example}
Let $\mathbb{C}1$ be the one dimensinal operator space.  Then
for $\alpha =[\alpha_{ij}] \in \mathbb{M}_n(\mathbb{C}1)$,
we have
$$
\mathcal{W}_{\text{max}}(\alpha) = w(\alpha).
$$
To see this, since
$
\mathcal{W}_{\text{max}}(\alpha)  = w([\alpha_{ij}z]) 
$
for some $z \in \mathbb{B}(\mathcal{K})$ with $\|z\|=1$, and
$\alpha$ double commutes with 
$\left[
\begin{array}{ccc}
z & & \\
 &\ddots & \\
 & & z \\
\end{array}
\right]
$,
we have  $\mathcal{W}_{\text{max}}(\alpha) \le w(\alpha).$
This and the maximality of $\mathcal{W}_{\text{max}}$ imply that
$$ 
w(\alpha)=
\inf \{\frac{1}{2} \| \beta \beta^* + \gamma^*\gamma \| \ | \ 
\alpha = \beta y \gamma , \|y\| =1, \  \beta, y, \gamma \in \mathbb{M}_n(\mathbb{C})\}.
$$
We note that the above equality for $w(\alpha)$ gives a simple proof of 
the Ando's Theorem in \cite{ando}, in case 
$\dim{\mathcal{H}} < \infty$.
\end{example}

\begin{example}
Let $X, Y$ be operator spaces in $\mathbb{B}(\mathcal{H})$.  For $x \in \mathbb{M}_{n, r}(X)$ and
$y \in \mathbb{M}_{r,n}(Y)$, we denote by $x \odot y$ the element 
$[\sum_{k=1}^r x_{ik} \otimes y_{kj}] \in  \mathbb{M}_{n}(X \otimes Y)$. We note that each
element $u \in \mathbb{M}_n(X \otimes Y)$ has a form $x \odot y$ for some $x \in \mathbb{M}_{n, r}(X)$,
$y \in \mathbb{M}_{r,n}(Y)$ and $r \in \mathbb{N}$.

\begin{flushleft}
{\bf{(a)}}
\end{flushleft}

We define
$$
\|u\|_{wh} =
\inf \{\frac{1}{2} \|xx^* + y^*y \| \ \mid \ u =x \odot y \in \mathbb{M}_n(X\otimes Y) \}
$$
for $u \in \mathbb{M}_n(X\otimes Y)$ (c.f. \cite{itohnagisa2}).
Then it is not hard to verify that $\| \ \|_{wh}$ satisfies the conditions 
WI and ${\text W}{\text{I{\hspace{-1pt}I}}}$. Moreover $\| \ \|_{wh}$ is a numerical
radius norm affiliated with $(X \otimes_h Y, \| \ \|_h)$, where $X \otimes_h Y$ is the
Haagerup tensor product operator space with the Haagerup norm $\| \ \|_h$, i.e.
 $$
\|u\|_{h} =
\inf \{ \|x\| \| y \| \ \mid \ u =x \odot y \in \mathbb{M}_n(X\otimes Y) \}.
$$

To see (OW),  given $u = x \odot y \in \mathbb{M}_n(X \otimes Y)$, we may assume that
$\| x \| = \| y \| $. Since
\begin{align*}
2 \left\| \left[
\begin{array}{cc}
0 & u \\
0 & 0 \\
\end{array} \right]
\right\|_{wh} 
&=
2\left\| \left[
\begin{array}{cc}
x & 0 \\
0 & 0 \\
\end{array}
 \right] \odot
 \left[
\begin{array}{cc}
0 & y \\
0 & 0 \\
\end{array}
 \right]
\right\|_{wh} \\
& \le 
\left\| \left[
\begin{array}{cc}
x & 0 \\
0 & 0 \\
\end{array}
 \right] 
\left[
\begin{array}{cc}
x & 0 \\
0 & 0 \\
\end{array}
 \right] ^*
+ 
\left[
\begin{array}{cc}
0 & y \\
0 & 0 \\
\end{array}
 \right]^*
 \left[
\begin{array}{cc}
0 & y \\
0 & 0 \\
\end{array}
 \right]
\right\| \\
& \le
\text{max} \{ \|x\|^2,  \|y \|^2 \} \\
& = \|x\| \|y\|, \\
\end{align*} 
we have $2 \left\| \left[
\begin{array}{cc}
0 & u \\
0 & 0 \\
\end{array} \right]
\right\|_{wh} \le \|u \|_h$. 

To see the other inequality, 
given $\varepsilon >0$. Since
$$
\left[\begin{array}{cc}
1 & 0\\
0 & 0 \\
\end{array}\right]
\left[\begin{array}{cc}
0 & u\\
0 & 0 \\
\end{array}\right]
\left[\begin{array}{cc}
0 & 0\\
0 & 1 \\
\end{array}\right]
=
\left[\begin{array}{cc}
0 & u \\
0 & 0 \\
\end{array}\right],
$$  
there exist $x_1, x_2, y_1, y_2$ such that
$
\left[\begin{array}{cc}
0 & u \\
0 & 0 \\
\end{array}\right]
=\left[\begin{array}{cc}
x_1 & x_2 \\
0 & 0 \\
\end{array}\right] \odot
\left[\begin{array}{cc}
0 & y_1 \\
0 & y_2 \\
\end{array}\right]
$.
Setting $x'=[x_1, x_2], 
y'=\left[\begin{array}{c}
y_1 \\
y_2\\
\end{array}\right]
$ with $\|x'\|=\|y'\|$, 
we rewrite 
$
\left[\begin{array}{cc}
0 & u \\
0 & 0 \\
\end{array}\right]
=\left[\begin{array}{cc}
x' & 0 \\
0 & 0 \\
\end{array}\right] \odot
\left[\begin{array}{cc}
0 & y' \\
0 & 0 \\
\end{array}\right]
$.
Thus we may assume that 
$u = x' \odot y'$ with $\max \{ \|x'{x'}^*\|, \|{y'}^*y' \| \} = \|x'\| \|y'\|$
and
$$
2\left\|
\left[\begin{array}{cc}
0 & u \\
0 & 0 \\
\end{array}\right] \right\|_{wh} + \varepsilon >
\left\| \left[\begin{array}{cc}
x' & 0 \\
0 & 0 \\
\end{array}\right]
\left[\begin{array}{cc}
x' & 0 \\
0 & 0 \\
\end{array}\right]^*
+
\left[\begin{array}{cc}
0 & y' \\
0 & 0 \\
\end{array}\right]^*
\left[\begin{array}{cc}
0 & y' \\
0 & 0 \\
\end{array}\right]
\right\|.
$$
Hence we obtain
$
2\left\|
\left[\begin{array}{cc}
0 & u \\
0 & 0 \\
\end{array}\right] \right\|_{wh}  \ge \|u \|_h.
$

\begin{flushleft}
{\bf{(b)}}
\end{flushleft}

We let denote  $X^\dagger = \{ x^* \in \mathbb{B}(\mathcal{H}) \ | \ x \in X \}$  
and also define a norm 
$\| \  \|_{wcb} $ on $X \otimes X^\dagger$ by
$$
\|u\|_{wcb} =
\inf \{\frac{1}{2} \|a\| \|x \|^2 \ \mid \ u =x a \odot x^* 
\in \mathbb{M}_n(X\otimes X^\dagger), 
x \in \mathbb{M}_{n, r}(X), a \in \mathbb{M}_r(\mathbb{C}), r \in \mathbb{N} \}
$$
for $u \in \mathbb{M}_n(X\otimes X^\dagger)$ 
(c.f. \cite{suen2}, \cite{itohnagisa1}).

It is easy to see that $\| \ \|_{wcb}$ also satisfies WI and 
${\text W}{\text{I{\hspace{-1pt}I}}}$. Since $\| \ \|_{wh}$ has another form \cite{itohnagisa2}
on $X \otimes X^\dagger$ as
$$
\|u\|_{wh} =
\inf \{ w(a) \|x \|^2 \ \mid \ u =x a \odot x^* 
\in \mathbb{M}_n(X\otimes X^\dagger), 
x \in \mathbb{M}_{n, r}(X), a \in \mathbb{M}_r(\mathbb{C}), r \in \mathbb{N} \},
$$
we have 
$$
\frac{1}{2} \|u \|_h \le \| u \|_{wcb} \le \| u \|_{wh} \le \| u \|_h, \quad
 u \in \mathbb{M}_n(X \otimes X^\dagger).
$$
Thus it turns out from Corollary 3.3 
that both $\| \ \|_{wh}$ and $\| \ \|_{wcb}$ are numerical radius norms
affiliated with the operator space $X \otimes_h X^\dagger$ with the Haagerup norm
$\| \ \|_h$.

\end{example}

We denote by $\mathcal{W}(X)$ the numerical radius operator space together with a
numerical radius norm $\mathcal{W}$ affiliated with an operator space $(X, \mathcal{O}_n)$. We call
$\mathcal{W}(X)$ a numerical radius operator space affiliated with $(X, \mathcal{O}_n)$.
Let $X, Y$ be operator spaces. It is clear that if $\varphi : X \longrightarrow Y$ 
is completely bounded,
then $\varphi :\mathcal{W}_{(1)}(X) \longrightarrow \mathcal{W}_{(2)}(Y)$ 
is $\mathcal{W}$-completely
bounded.

\begin{lem}Let $X, Y$ be operator spaces and $\mathcal{W}(X)$ a numerical 
radius operator space affiliated with $X$.
If $\varphi : X \longrightarrow Y$ is completely bounded, then
$\mathcal{W}
( \varphi : \mathcal{W}(X) \longrightarrow \mathcal{W}_{\text{min}}(Y))_{cb} 
\le \mathcal{O}(\varphi)_{cb}$.
\end{lem}

\begin{proof}
It follows from
\begin{align*}
\mathcal{W}(\varphi)_{cb} 
& = \sup \{ \frac{1}{2} \mathcal{O}(\varphi_n(x)) \ | 
        \ \mathcal{W}(x) \le 1,  x \in \mathbb{M}_n(X), n \in \mathbb{N} \}\\
& \le \sup \{ \frac{1}{2} \mathcal{O}(\varphi_n(x) )\ | 
          \frac{1}{2}\mathcal{O}(x) \le 1, x \in \mathbb{M}_n(X), n \in \mathbb{N} \}\\
&= \mathcal{O}(\varphi)_{cb}.
\end{align*}
\end{proof}

\begin{lem}Let $X, Y$ be operator spaces and $\mathcal{W}(Y)$ a numerical 
radius operator space affiliated with $Y$.
If $\varphi : X \longrightarrow Y$ is completely bounded, then
$\mathcal{W}
( \varphi : \mathcal{W}_{\text{max}}(X) \longrightarrow \mathcal{W}(Y))_{cb} 
\le \mathcal{O}(\varphi)_{cb}$.
\end{lem}
\begin{proof}
Assume that $\mathcal{O}(\varphi)_{cb} \le 1.$ 
Let $x \in \mathbb{M}_n(\mathcal{W}_{\text{max}}(X))$ with $\mathcal{W}_{\text{max}}(x) \le 1$.
Since $\mathcal{W}(Y)$ has a $\mathcal{W}$-complete isometry 
$\varPhi : \mathcal{W}(Y) \longrightarrow w(\mathbb{B}(\mathcal{H}))$,
we have $\mathcal{W}(\varphi_n(x)) =w(\varPhi_n \circ \varphi_n(x))$. 
We note that
$w(\varPhi_n \circ \varphi_n(x)) \le \mathcal{W}_{\text{max}}(x)$, since 
$\mathbb{M}_n(\mathcal{W}_{\text{max}}(X/ \ker ( \varPhi \circ \varphi)) ) \ni \tilde{x} \longmapsto 
\varPhi_n \circ \varphi_n(x) \in \mathbb{M}_n(w(\mathbb{B}(\mathcal{H})))$ is isometric.
Hence we have $\mathcal{W}(\varphi_n(x)) \le 1$.  
\end{proof}

\begin{lem}Let $X, Y$ be operator spaces and $\mathcal{W}(X), \mathcal{W}(Y)$ numerical 
radius operator spaces affiliated with $X, Y$.
If $\varphi : X \longrightarrow Y$ is completely bounded, then
$\mathcal{O}(\varphi)_{cb} \le \mathcal{W}
( \varphi : \mathcal{W}(X) \longrightarrow \mathcal{W}(Y))_{cb} $.
\end{lem}

\begin{proof}

It follows from
\begin{align*}
 & \quad \mathcal{O}(\varphi)_{cb} \\
& = \sup \{ 2 \mathcal{W}\left(
\left[\begin{array}{cc}
0 & \varphi_n(x)\\
0 & 0 \\ 
\end{array}\right]
\right) \ | 
        \ 2 \mathcal{W}\left(
\left[\begin{array}{cc}
0 & x \\
0 & 0 \\ 
\end{array}\right]\right)
 \le 1,  x \in \mathbb{M}_n(X), n \in \mathbb{N} \}\\
& \le \sup \left\{  \mathcal{W}(\varphi_{2n}(y) )\ | 
          \mathcal{W}(y) \le 1, y \in \mathbb{M}_{2n}(X), n \in \mathbb{N} \right\}\\
&= \mathcal{W}(\varphi)_{cb}.
\end{align*}
\end{proof}

We let $\mathbb{O}$ denote the category of operator spaces, in which the objects are
the operator spaces and the morphisms are the completely bounded maps. We also
let $\mathbb{W}$ denote  the category of numeical radius operator spaces with 
the morphisms being the $\mathcal{W}$-completely bounded maps.
We have already obtained a functor $\mathcal{O} : \mathbb{W} \longrightarrow \mathbb{O}$
such that  $\mathcal{O}(X) = 2 \mathcal{W}
\left(
\begin{array}{cc}
0 & X \\
0 & 0 \\ 
\end{array}\right)$ symbolically. We have also found functors 
$\mathcal{W} : \mathbb{O} \longrightarrow \mathbb{W}$ which satisfy
$\mathcal{O}\circ \mathcal{W}(X) = X$ for each operator space $X$.
Combining the above Lemmas, $\mathcal{W}_{\text{max}}$ and $\mathcal{W}_{\text{min}}$
can be seen as the functors which embed $\mathbb{O}$ into $\mathbb{W}$ strictly.

\begin{thm}
Let $X, Y$ be operator spaces. If $\varphi : X \longrightarrow Y$ is a linear map, then
\begin{enumerate}
\item[(1)]$ \mathcal{W}( \varphi : \mathcal{W}_{\text{max}}(X) 
\longrightarrow \mathcal{W}_{\text{max}}(Y))_{cb} =
\mathcal{O}(\varphi : X \longrightarrow Y )_{cb}$,
\item[(2)] $\mathcal{W}( \varphi : \mathcal{W}_{\text{min}}(X) 
\longrightarrow \mathcal{W}_{\text{min}}(Y))_{cb} =
\mathcal{O}(\varphi : X \longrightarrow Y )_{cb}$.
\end{enumerate}

\end{thm}

\begin{proof}
It is clear from Lemma 3.6, Lemma 3.7 and Lemma 3.8.

\end{proof}

\end{document}